\documentclass[12pt,a4paper]{article}  
\usepackage{t1enc}
\usepackage[latin1]{inputenc}
\usepackage[english]{babel}
\usepackage{graphicx}
\usepackage{amsmath}
\usepackage{amssymb}
\usepackage{amsthm}
\usepackage{pgfplots}
\usepackage{tikz-3dplot}
\usepackage{tikz}
\usepackage{ifthen}

\usetikzlibrary{cd, arrows, matrix, intersections, math, calc}
\font\bBB=msbm10

\def\bBR{\mbox{\bBB R}}

\def\bBZ{\mbox{\bBB Z}}

\def\R3{\bBR ^3}

\def\Zz{\bBZ_2}

\begin{document}
\title {Topology of $SO(3)$ for Kids}
\author{Orlin Stoytchev\thanks{American University in Bulgaria, 2700 Blagoevgrad, Bulgaria} }
\date{}
\maketitle
\abstract {We present a simple visual description of the topology of the space of three-dimensional rotations, requiring just intuition, imagination and no advanced math.  }
\section{Introduction}
By a {\it rotation in 3-space}, or a {\it three-dimensional rotation} we may understand changing the orientation of an object. While rotating an object by 360$^{\circ}$ about any axis brings the object back to its initial orientation, the motion we have performed "leaves a trace". If we attach a band (belt) to the object and fix its other end, after the  360$^{\circ}$ rotation the band will be twisted and cannot be untwisted without further rotating the object. Our world, however, has a fascinating property - if we rotated the object by 720$^{\circ}$, the band can be untwisted by flipping it around the object. What we have described  is the so-called "Dirac's belt trick" (see, e.g., \cite{Egan}, \cite{Palais} for nice Java applets and explanations). This is related to the fact that a continuous motion which is a rotation from 0$^{\circ}$ to 360$^{\circ}$ is not continuously deformable to no motion at all (i.e., trivial motion), while a rotation from 0$^{\circ}$ to 720$^{\circ}$ is. We say that the first motion is topologically nontrivial (in some sense, to be discussed later), while the second is topologically trivial.\par
The set of all possible rotations in three-dimensional space is denoted by $SO(3)$. It is a topological space since there is an obvious natural notion of what it means for two points to be close to each other. Moreover, it is a closed smooth three-dimensional manifold, i.e., a three-dimensional analog of a closed smooth surface. In addition $SO(3)$ is a group - you can compose rotations (multiplication) and each rotation has an inverse. This group has a concrete {\it representation} as the group of $3\times 3$ orthogonal matrices with determinant 1, hence the name ``SO(3)'' = ``special orthogonal  matrices of size 3''. The group structure and the structure of smooth manifold are compatible with each other, so this is what is called a Lie group.\par
 The ``fascinating property'' referred to above is linked to a topological property of $SO(3)$. (A topological property is one which is preserved under deformations of the space. No cutting and gluing is allowed. If we cut our space at some intermediate step, we should later glue back the points along the cut.) Taking an object and continuously changing its orientation, starting and ending at the initial orientation (a {\it full rotation}) corresponds to a closed path in SO(3), starting and ending at the identity. It turns out that there are two classes of full rotations: topologically trivial, i.e., deformable to the trivial motion, and topologically nontrivial. Every nontrivial full rotation can be deformed to any other nontrivial full rotation. Two consecutive nontrivial full rotations produce a trivial one.\par
For any topological space, one can consider the set of closed paths starting and ending at some fixed point, called {\it base-point}. Two closed paths that can be continuously deformed to each other, keeping the base-point fixed, are called \textit{homotopic}. One can multiply closed paths by {\it concatenation}, i.e., take the path obtained (after appropriate reparametrization) by traveling along the first and then along the second. There is also an inverse for each path--the path traveled in reverse direction. These operations turn the set of homotopy classes of closed paths (with a given base-point) into a group and it is an important topological invariant of any topological space. It was introduced by Poincar\'{e} and is called the {\it first homotopy group} or the \textit{fundamental group} of the space, denoted by $\pi_1$. Thus, the property of $SO(3)$ stated above is written concisely as $\pi_1(SO(3)) \cong \bBZ/2\bBZ\equiv \Zz$. In other words, $\pi_1(SO(3))$ is isomorphic to the group of two elements $\{1,-1\}$ with the usual multiplication.\par
This specific topological property of $SO(3)$ plays a fundamental role in our physical world. Because there are two homotopy classes of closed paths in $SO(3)$ there exist precisely two principally different types of elementary particles: those with integer spin (e.g., the photon), which turn out to be bosons, and those with half-integer spin (e.g., the electron), which are fermions. The difference can be traced to the fact that the complex (possibly multicomponent) wave function determining the quantum state of a boson is left unchanged by a rotation by $2\pi$ of the coordinate system while the same transformation multiplies the wave function of a fermion by $-1$. This is possible since only the modulus of the wave function has a direct physical meaning, so measurable quantities are left invariant under a full rotation by $2\pi$. However, as discovered by Pauli and Dirac, one needs to use wave functions having this (unexpected) transformation property for the correct description of particles with half-integer spin, such as the electron. The careful analysis showed (\cite{Wig,Barg}) that the wave function has to transform properly only under transformations which are in a small neighborhood of the identity. A "large" transformation such as a rotation by $2\pi$ can be obtained as a product of small transformations, but the transformed wave function need not come back to itself--there may be a complex phase multiplying it. From continuity requirements it follows that if one takes a closed path in $SO(3)$ which is contractible, the end-point wave function must coincide with the initial one. Therefore, a rotation by $4\pi$ should bring back the wave function to its initial value and so the phase factor corresponding to a $2\pi$-rotation can only be $-1$. What we have just described is the idea of the so-called {\it projective} representations of a Lie group, which we have to use in quantum physics. As we see on the example of $SO(3)$, they exist because the latter is not simply-connected, i.e., $\pi_1(SO(3))$ is not trivial. Projective representations of a non-simply-connected Lie group are in fact representations of its covering groups. In the case of $SO(3)$ this is the group $SU(2)$ of 
$2\times 2$ unitary matrices with unit determinant. Topologically, this is the three-dimensional sphere $S^3$; it is a double cover of $SO(3)$ and the two groups are locally isomorphic and have the same Lie algebra.
\par 
A refinement of "Dirac's belt trick", which may be called ``braid trick'' can be used to actually show that $\pi_1(SO(3)) \cong \bBZ/2\bBZ\equiv \Zz$. The idea is fairly simple and is based on the following experiment: attach the ends of three strands to a ball, attach their other ends to the desk, perform an arbitrary number of full rotations of the ball to obtain a plaited braid. Then try to unplait it without further rotating the ball. As expected, braids that correspond to contractible paths in $SO(3)$ are trivial, while those corresponding to noncontractible paths form a single nontrivial class. In \cite{Sto} a one-to-one correspondence was constructed between homotopy classes of closed paths in $SO(3)$ and a certain factor group of the group $P_3$ of pure braids with three strands. This factor group turns out to be isomorphic to $\bBZ/2\bBZ$.\par
The standard way of showing that $\pi _1(SO(3))\cong \bBZ _2$ is to prove that $SO(3)$ is topologically equivalent ({\it homeomorphic}) to $S^3/i$, where $S^3$ denotes the tree-dimensional sphere and
$i$ means identifying diametrically opposite points of $S^3$. In topology the space $S^3/i$ is also called the {\it real three-dimensional projective space}, denoted by $\bBR\text P^3$.  $\bBR\text P^3$ is formally defined as the space of all lines in $\bBR^4$. Notice that lines in $\bBR^4$ are in one-to-one correspondence with pairs of opposite points on the sphere $S^3$.  Once the latter is known, it is quite easy to see that paths starting from one point of $S^3$ and ending at the opposite point will be closed paths in $S^3/i$, which are not contractible, but the composition of two such paths gives a contractible path. \par
In order to find the topological structure of $SO(3)$ one normally uses Lie group and Lie algebra theory. Namely, it is shown that the group $SU(2)$ of unitary $2\times 2$ matrices with determinant 1 is homeomorphic to $S^3$ and that there is a 2--1homomorphism $SU(2)\to SO(3)$, which is a local isomorphism, and which sends diametrically opposite points in $SU(2)$ to the same point in $SO(3)$. The interested reader may find a nice treatment of this and many other aspects of three-dimensional rotations in \cite{Egan2}. \par
A more direct and geometric way to exhibit the topology of $SO(3)$ uses the fact that any rotation is a rotation about some fixed axis. This is neither intuitive nor so easy to prove and is known as Euler's Rotation Theorem \cite{Eu}, \cite{Wiki}. Nowadays a standard approach to prove it is to use linear algebra and that rotations are represented by $3\times 3$ matrices and any such matrix must have one (real) eigenvalue and a corresponding eigenvector. Then one shows easily that this element of $SO(3)$ must be a rotation about this eigenvector. In other words  we need to specify the angle of rotation and the orientation of that axis in $\bBR^3$ with the necessary identification with rotation about the opposite axis at the opposite angle (see also \cite{Peng}). \par
Our visual derivation of the topology of $SO(3)$ uses a different paramerization of rotations and does not rely on results from Lie group theory, linear algebra, or Euler's theorem, and can be accessible to anyone with enough geometric intuition.\par

\section{Visualizing the Space $SO(3)$}
\begin{figure} [h!]
\centering
\includegraphics[width=70 mm]{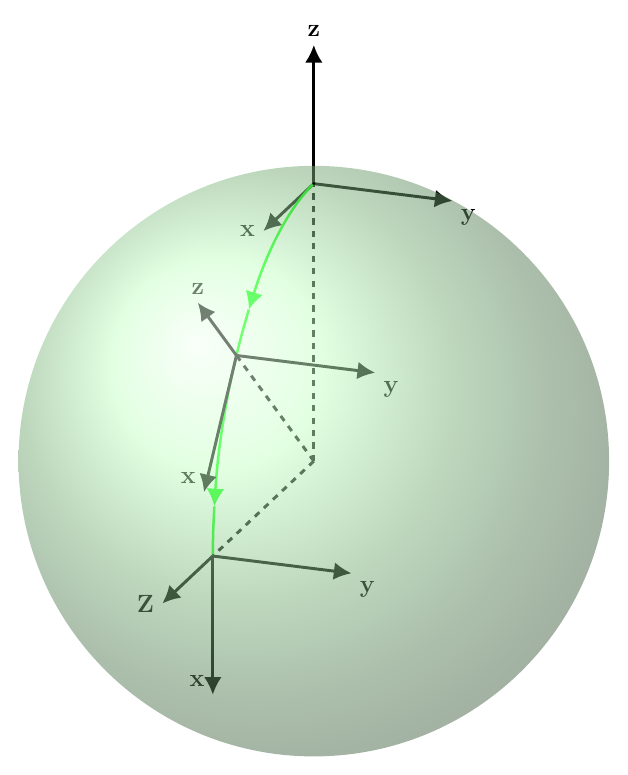}
\caption{Rotating the z-axis to the Z-axis}
\label{paral}
\end{figure}

It is standard to consider any rotation as a rotation of a rectangular right-handed coordinate system. Thus we have the original coordinate axes $\text{x,\,y,\,z}$ and the rotated ones $\text{X,\,Y,\,Z}$. We can first rotate $\text{z}$ to $\text{Z}$ along a shortest geodesic on a unit sphere (Figure \ref{paral}).
Geodesics are locally shortest curves between points. On the sphere with its standard coordinate grid geodesics are all the meridians and the equator. Initially the z-axis points to the north pole. As we move the z-axis along a meridian we will agree to perform parallel transport of the x- and y-axes. For convenience we have shifted these two axes  to a distance 1 from the origin, so they are both tangent to the sphere. Parallel transport of a tangent vector along a geodesic means moving the vector along the geodesic while keeping the angle between it and the geodesic fixed. Let's specify that, e.g.,  the x-axis points along the $0^{\circ}$ meridian (the "Greenwich meridian"), which is the one drawn in the picture. The prescribed recipe works in all cases except when the Z-axis is in the direction of the south pole, so we exclude this case for the time being . \par
\begin{figure} [h!]
\centering
\includegraphics[width=70 mm]{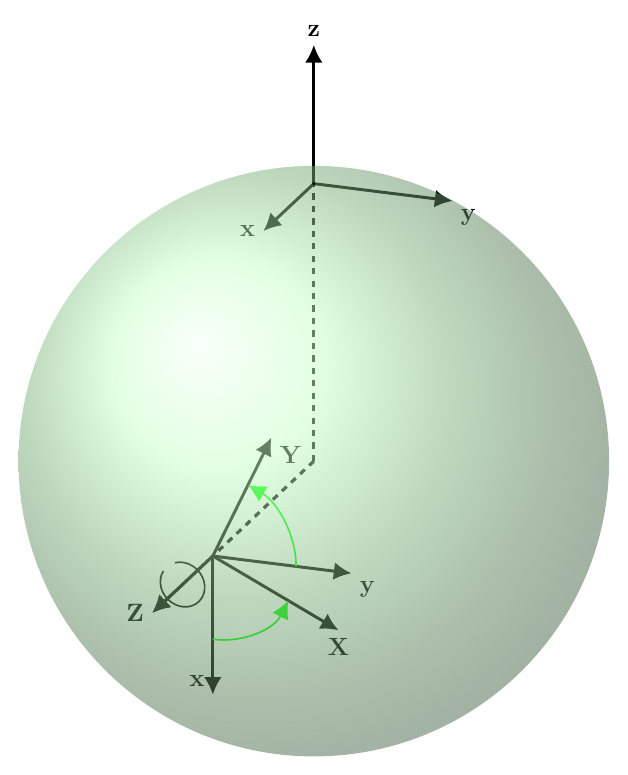}
\caption{Rotating about the Z-axis}
\label{rotat}
\end{figure}

Next we rotate about the Z-axis until the intermediate x- and y-axes match X and Y (Figure \ref{rotat}).
All such rotations are specified by a point on the punctured unit sphere with the south pole removed and a point on a unit circle. The punctured sphere is topologically the same as the open disk. Thus the space of all rotations except those where Z is opposite to z is topologically equivalent to the interior of a solid torus.\par
\begin{figure} [h!]
\centering
\includegraphics[width=70 mm]{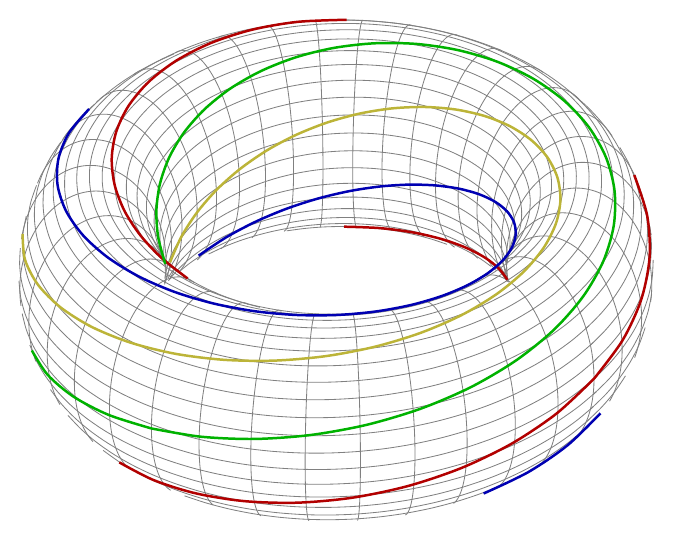}
\caption{$SO(3)\cong$ solid torus with  points on the boundary identified along winding curves}
\label{torus}
\end{figure}

When the Z-axis passes through the south pole we can choose any meridian to slide along. It is easy to see that sliding along the $90^{\circ}$-meridian will produce orientation of the x- and y-axes which is rotated by $180^{\circ}$ relative to their orientation if we slide along the ``Greenwich'' meridian. Further, sliding along the  $180^{\circ}$-degree meridian produces x- and y-axes rotated by $360^{\circ}$ relative to the latter. As we continue increasing the angle of the meridian that we use, the x- and y-axes perform one more full rotation. This means that we have to identify points on the boundary of the solid torus along curves which for one turn along the small circle perform two turns along the large circle. In Figure \ref{torus} all points belonging to  a curve on the surface parallel to the colored ones should be identified.

\begin{figure} [h!]
\centering
\includegraphics[width=70 mm]{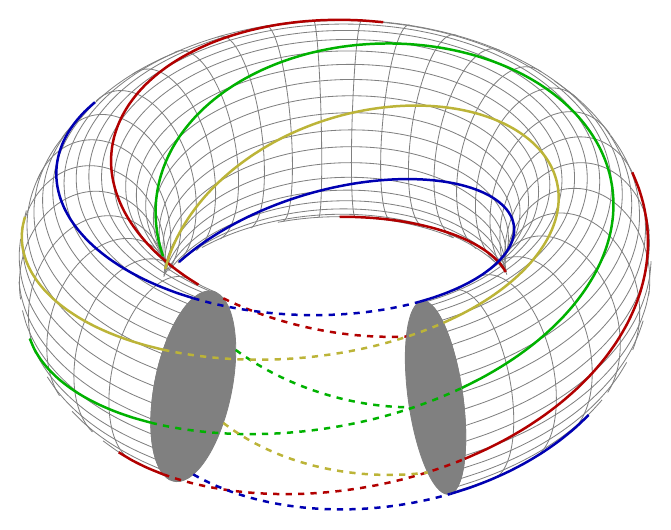}
\caption{Cutting the torus along a disk}
\label{cuttorus}
\end{figure}

It turns out that the solid torus with the specified identification of points on the boundary is in fact topologically equivalent to $\bBR\text P^3$. To see this we first (temporarily) cut the torus along a disk, as shown on Figure \ref{cuttorus}.

\begin{figure} [h!]
\centering
\includegraphics[width=70 mm]{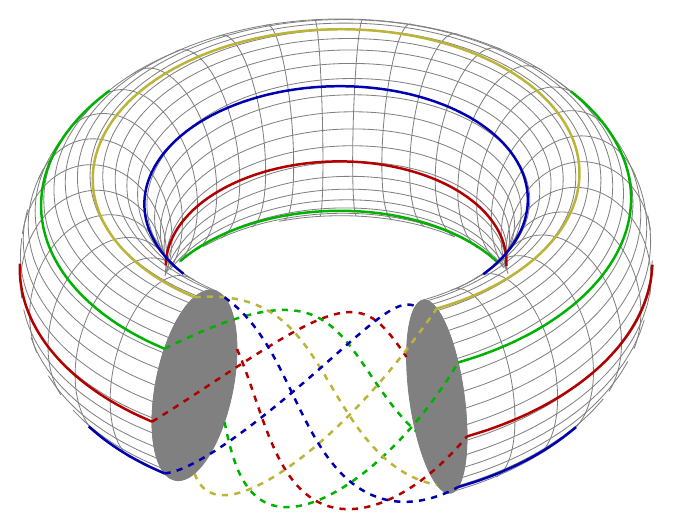}
\caption{Twisting the torus}
\label{twistedtorus}
\end{figure}

\begin{figure} [h!]
\centering
\includegraphics[width=25 mm]{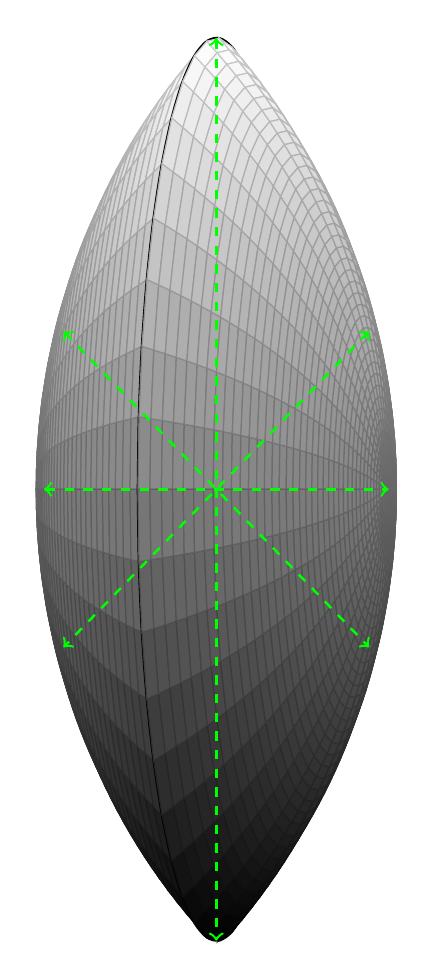}
\caption{Shrinking the boundary of the cut torus}
\label{lens}
\end{figure}

Next we twist the solid torus by $180^{\circ}$ so that the colored curves run straight along the surface. The dashed curves between the two sides of the cut indicate how we should identify back points along the cut later (Fig. \ref{twistedtorus}). Note that each point on one side of the cut should be identified with a point on the other side, rotated by $180^{\circ}$.
Then we shrink the boundary of the cut torus along curves parallel to the colored ones to obtain a lens-shaped space. The green dashed lines on Figure \ref{lens} indicate how we have to identify opposite points. Now it is easy to see that this is in fact topologically a ball with diametrically opposite points on the boundary identified.\par
It is well-known that this is topologically the same as  $\bBR\text P^3$. Since our aim is to make the exposition as self-contained as possible, we sketch the explanation. The unit three-dimensional sphere is $S^3=\{(x_1,x_2,x_3,x_4)\in \bBR^4\ |\ { x_1}^2+{x_2}^2+{x_3}^2+{x_4}^2=1\}$. The lower hemisphere $S^{3}_-$ can be defined as the subset of $S^3$ for which $x_4\le 0$. For each point on the lower hemisphere there is a diametrically opposite point on the upper hemisphere which is identified with the former. So we can remove the open upper hemisphere and we get the lower hemisphere with diametrically opposite points on the boundary identified. Notice that the boundary of $S^{3}_-$ consists of those points on $S^3$ for which $x_4=0$, i.e., this is the unit two-dimensional sphere $S^2$. Finally we deform (collapse) the lower three-dimensional hemisphere by setting $x_4$ to zero to get a three-dimensional ``disk'', i.e., a usual three-dimensional ball. \par
\begin{figure}[h!]
\centering
\includegraphics[width=40 mm]{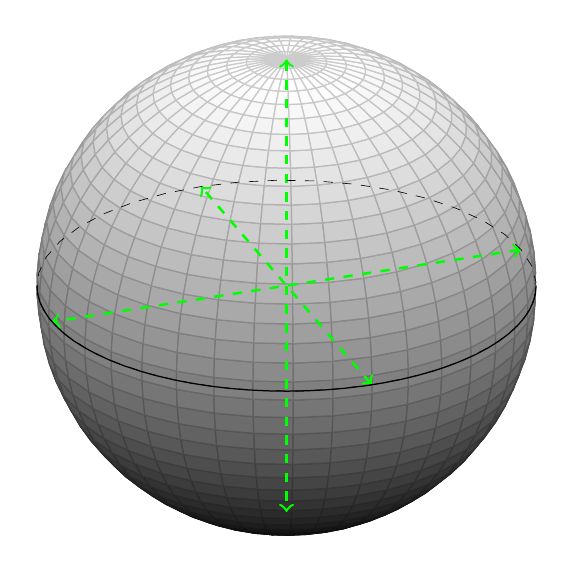}
\includegraphics{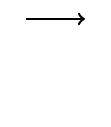}
\includegraphics[width=40 mm]{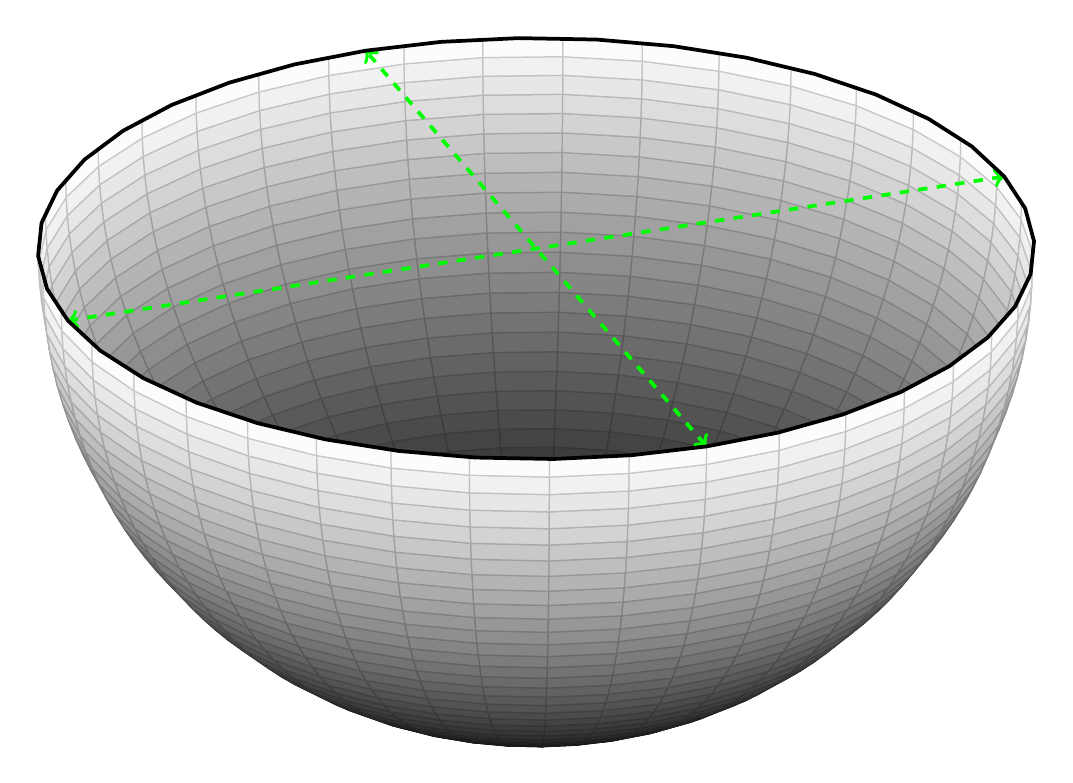}
\includegraphics{arrow.pdf}
\includegraphics[width=40 mm]{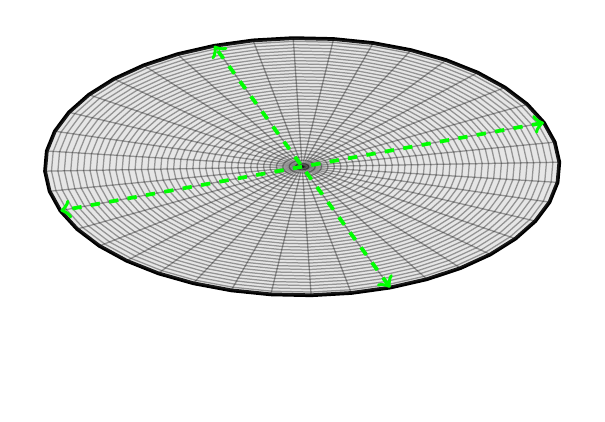}
\caption{$\bBR{\text P}^2 \cong$ disk with diametrically opposite points identified}
\label{defform}
\end{figure}

We illustrate the process with pictures in one dimension less - the two-dimensional sphere with diametrically opposite points identified is topologically the same as the disk with opposite points on the boundary identified (Figure \ref{defform}).
\section{Contractible and Non-Contractible Closed Paths in $SO(3)$}
\begin{figure}[h!]
\centering
\includegraphics[width=40mm]{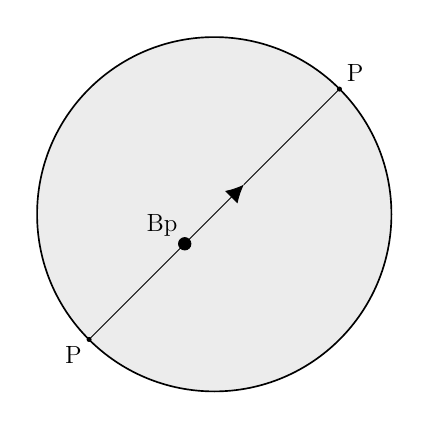}
\caption{A non-contractible path}
\label{paths1}
\end{figure}

\begin{figure}[h!]
\centering
\includegraphics[width=40 mm]{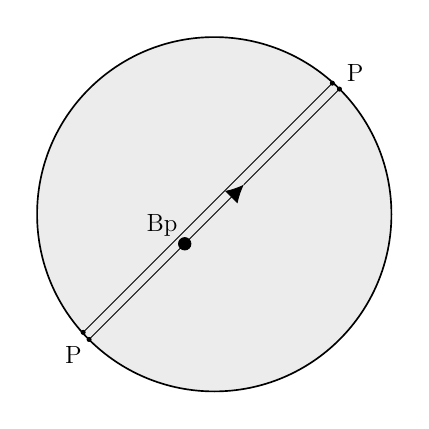}
\includegraphics{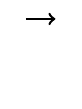}
\includegraphics[width=40 mm]{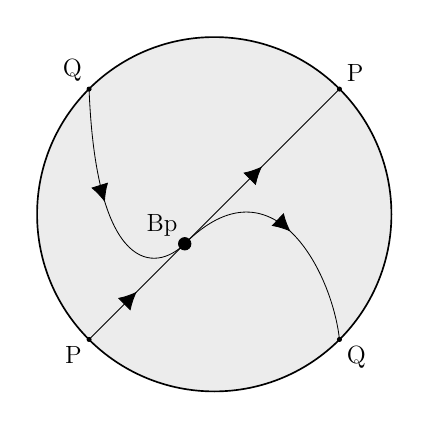}
\includegraphics{arrow2.pdf}
\includegraphics[width=40 mm]{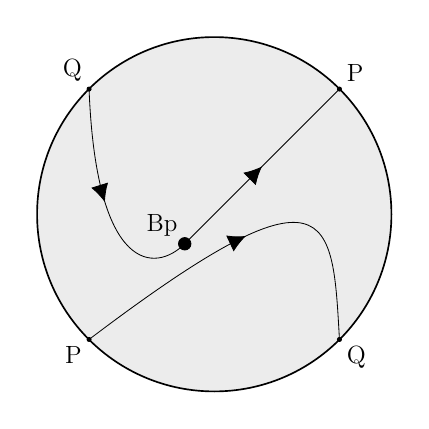}
\includegraphics{arrow2.pdf}
\includegraphics[width=40 mm]{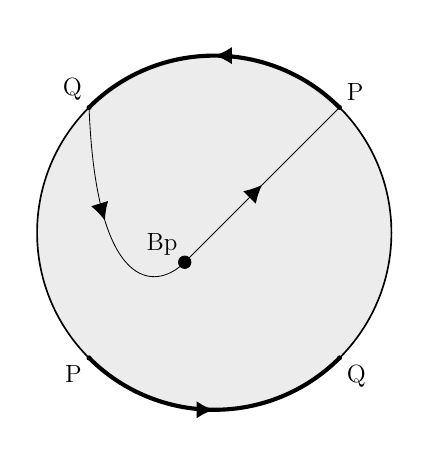}
\includegraphics{arrow2.pdf}
\includegraphics[width=40 mm]{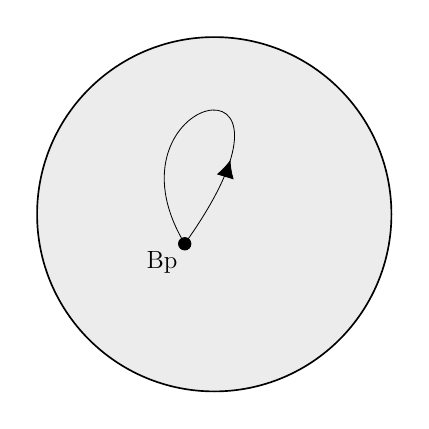}
\includegraphics{arrow2.pdf}
\includegraphics[width=40 mm]{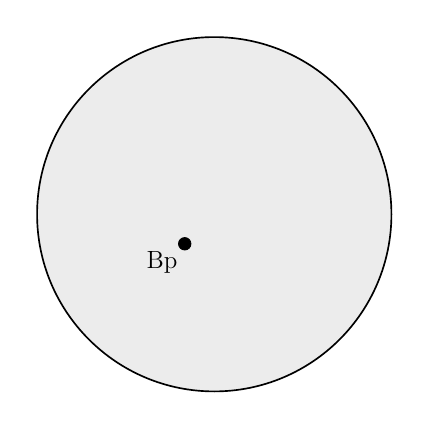}
\caption{A contractible path}
\label{paths2}
\end{figure}

Having identified the topology of $SO(3)$ it is easy to understand the nature of contractible and non-contractible closed paths. For example a path as shown on Figure \ref{paths1} starting from a base-point Bp in the interior of the ball, moving along a diameter to the boundary point P and coming from the opposite side back to Bp is evidently not contractible.
However, traversing the same path twice we get a contractible path, as illustrated in Figure \ref{paths2}. Note that the middle point is again Bp but now we are allowed to move it.

\end{document}